\documentclass[11pt]{article}
\usepackage{epsfig}
\usepackage[utf8]{inputenc} 

\usepackage{rotating}
\usepackage{graphicx}
\usepackage{verbatim}

\usepackage{url}
\usepackage{color}
\usepackage{pstricks}
\usepackage{fancyvrb}
\usepackage{makeidx}

\VerbatimFootnotes

\begin{document}

\newcommand{\bm}[1]{\mbox{\boldmath$#1$}}

\def\mvec#1{{\bm{#1}}}   

\title{More lessons from the six box toy 
experiment\footnote{Addendum to {\tt arXiv:1612.05292}.} 
}

\author{G.~D'Agostini \\
Universit\`a ``La Sapienza'' and INFN, Roma, Italia \\
{\small (giulio.dagostini@roma1.infn.it,
 \url{http://www.roma1.infn.it/~dagos})}
}

\date{}

\maketitle

\begin{abstract}

Following a paper in which the fundamental aspects
of probabilistic inference were introduced by means of a toy
experiment, details of the analysis of simulated
long sequences of extractions are shown here.
In fact, the striking performance of probability-based
inference and forecasting, compared to those
obtained by simple `rules', 
might impress those practitioners who are usually 
underwhelmed by the philosophical foundation 
of the different methods.
The analysis of the sequences also shows how the 
smallness of the probability of what has actually been 
 observed, given the hypotheses of interest, 
is irrelevant for the purpose of inference.
\end{abstract}


\begin{flushright}
{\sl ``Grown-ups like numbers''}\\
(Saint-Exup\'ery's Little Prince)
\end{flushright}

\section{Introduction}
For years I have been  using a toy experiment
for introducing probabilistic reasoning. 
Irrespective of whether my audience has been of professional 
physicists and engineers,\footnote{See e.g. 
\url{http://indico.cern.ch/event/419045/}
and \url{http://www.roma1.infn.it/~dagos/prob+stat.html#cern05}; 
\url{http://2015.imtc.ieee-ims.org/content/tutorials} 
and \url{http://www.roma1.infn.it/~dagos/prob+stat.html#IEEEPisa}.
} 
high school 
students, teachers and general public,\footnote{See e.g. 
\url{http://www.lnf.infn.it/edu/openlabs/2016/conference.php} 
and
\url{http://orientamento.matfis.uniroma3.it/fisincittastorico.php#dagostini}
(in Italian).
} 
or even managers and senior officers attending a 
decision-making school,\footnote{See e.g. 
\url{http://www.pangeaformazione.it/en/training/decision-making-school.html}.} 
this toy model has always been an ``eye opener''. 
This is how it was defined by the editors of a special issue 
of the American Journal of
 Physics,\footnote{\url{http://stp.clarku.edu/ajp_contributors.html}.} 
in which this `experiment' was first published \cite{AJP}. 

A thorough description of the game, and what can be learned from it, 
is given in Ref. \cite{ME2016}.
In particular, in that paper I explain the reasons 
why I do not let anyone see the balls in the box.  
At most we can make simulated extractions, 
in which we can `almost' see the game from the 
God's perspective: we know the box composition
with certainty, and give a
superior smile at the algorithm that is trying 
to guess it.
`Almost', because we remain uncertain 
about the color of future extractions. 

Let us see then what we can learn from simulations. 
Firstly, in order to allow readers to reproduce the results,
details of the simulation are given in Section 2. 
For this reason, as examples of how to generate and analyze 
the sequences, R commands are given. 
In Sections 3-5, some sequences are analyzed in detail. Whenever possible,
the numbers obtained from the probability theory algorithm are compared
with those resulting from `simple rules'. But, as made clear in
Section 6, this is not always possible. 
Finally,  in Section 7, I emphasize the fact that most real life cases --
as random sequences of black and white balls extracted from a box are -- 
might have `astronomically' small probabilities of occurrence,
given the hypotheses of interest. 
But, nevertheless, the smallness of each  
conditional probability is irrelevant for the 
inference. Instead, what matters are their ratios 
and the relative prior beliefs of the different hypotheses. 

\section{Simulated sequences}
\begin{figure}
\centerline{\includegraphics[width=0.85\linewidth]{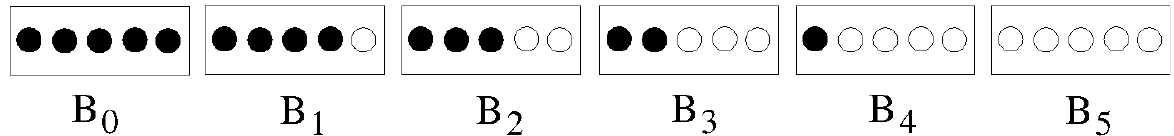}}
\mbox{}\vspace{-0.3cm}
\caption{\small \sf A sketch of the six boxes of the toy experiment.
The index refers to the number of white balls.}
\label{fig:sixbox}
\end{figure}
Thousand extractions from each 
of the  boxes $B_0$, $B_1$ and $B_2$
(Fig.~\ref{fig:sixbox}), since we can form an idea
of what happens from the others just by (anti-)symmetry. 
The R code to generate and analyze the sequences 
is based on that shown in Footnote 31
of \cite{ME2016}, but we report here also the 
inferential story as the  extractions go on. 
Moreover, for the benefit of the reader, 
who can then check the details of the
results the `seed' of the random generator is given,
equal to 20160715, for the date of the talk upon which
this paper is based (no special, fancy sequences have been
cherry-picked).

Here are the four lines of R code to make initializations 
and extractions (`0' for Black and `1' for
 White):\footnote{The csv files of the sequences analyzed here
can be downloaded from 
\url{http://www.roma1.infn.it/~dagos/prob+stat.html#addendum_ME}, 
and then loaded in a R session 
by e.g. \verb|seq.B0 = as.vector( read.csv("seq_B0.txt")$x )|.
}\\
\verb|N = 5; i = 0:N; pii = i/N; n = 1000; set.seed(20160715)|\\
\verb|seq.B0 = rbinom(n, 1, pii[1])|\\
\verb|seq.B1 = rbinom(n, 1, pii[2])|\\
\verb|seq.B2 = rbinom(n, 1, pii[3])|\\
Later, in order to have a feeling of the performances
of the method,\footnote{I you wish and it helps you, 
you might think of the `propensity' \cite{ME2016} of the algorithm 
to produce some numbers rather then others.}
we can split the sequences in `runs' of 100 extractions
and analyze them independently. 

\section{Box $B_0$}
\begin{figure}
\centerline{\includegraphics[width=0.9\linewidth]{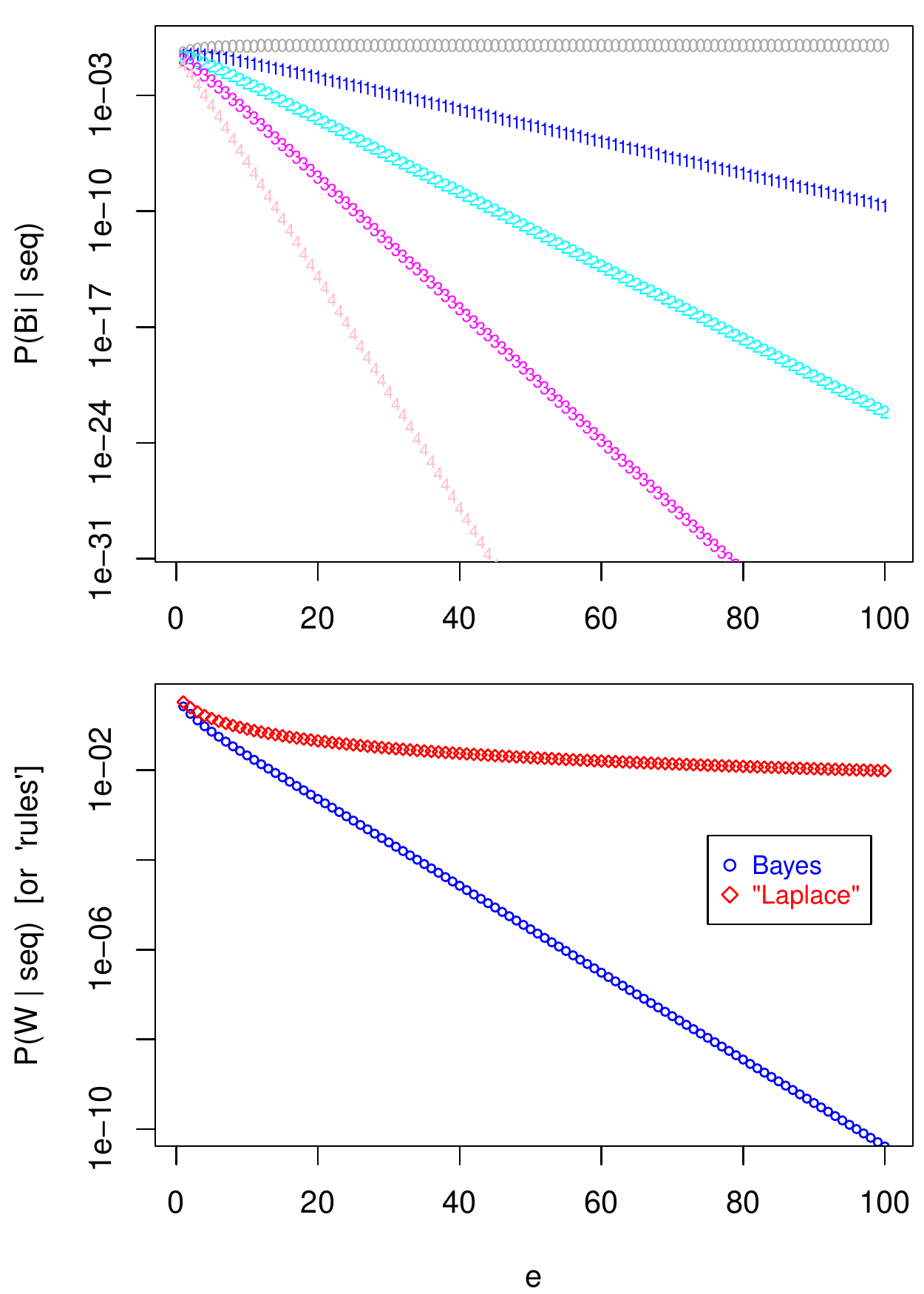}}
\caption{\small \sf Simulated extractions from $B_0$ (see text, 
especially for the meaning of "Laplace", that indeed stands for
"\underline{misused} Laplace {\em rule of succession}").}
\label{fig:Simulazione_B0}
\end{figure}

The runs from $B_0$ 
will obviously be all equal, since that composition can only 
produce 0's (Black), and therefore only the first run is shown
(Fig. \ref{fig:Simulazione_B0}). 
We see (upper plot) that as we continue 
to observe Black our confidence  that we have 
`picked' $B_0$ steady increases. $B_5$ is obviously ruled
out from the very beginning, while the probability
we {\em have to} assign
the compositions having both colors
 decreases exponentially
 with the increasing number of extractions 
(note the log scale in the ordinates).\footnote{The reason is quite simple,
with the approximation being
valid for `large' values of $n$ (and for $i\ne 5$):
\begin{eqnarray*}
P(B_i\,|\,n\mbox{B},I) &=& \frac{P(n\mbox{B}\,|\,B_i,I)\times 1/6}
    {1\times 1/6 + (4/5)^n\times 1/6 + \cdots} 
= \frac{(5-i)/5\times 1/6}
    {1\times 1/6 + (4/5)^n\times 1/6 + \cdots} \\
 &\approx& \left(\frac{5-i}{5}\right)^n = \left(\frac{5}{5-i}\right)^{-n}\,.   
\end{eqnarray*}
\label{fn:exponential_PBi},
 }
After 100 extractions the probability of $B_0$ differ from unity 
by about $10^{-10}$, essentially the probability of $B_1$, while 
all others are less probable by tens of orders of 
magnitude.\footnote{More precisely, executing the R command
of Footnote 31 of Ref. \cite{ME2016} with \verb|ri=0| we get \\
\verb|[1] 1.000000e+00 2.037036e-10 6.533186e-23 1.606938e-40 1.267651e-70|
} So, we are {\em practically certain} about $B_0$ --  but 
those who only like certainty have to remember
that {\em our only certainty is that $B_5$ is ruled out}.

The lower plot of Fig. \ref{fig:Simulazione_B0} shows, instead,
the probability of White in a next extraction, that
is $P(\mbox{W}\,|\,\mbox{seq})$ (blue circles). Its exponential
decrease results 
from the exponential decrease of 
$P(B_i\,|\,n\mbox{B},I)$, for $i>0$.\footnote{In fact, 
using the result of Footnote
\ref{fn:exponential_PBi}, 
\begin{eqnarray*}
P(\mbox{W}\,|\,n\mbox{B},I) &=& 
  P(\mbox{W}\,|\,B_1,I)\cdot P(B_1\,|\,n\mbox{B},I) + \cdots \\
   &\approx& \frac{1}{5}\times \left(\frac{5}{4}\right)^{-n} \,.
\end{eqnarray*}
With $n=100$ we obtain $4\times 10^{-11}$, or 1 in 25 billions
(remember that `who'  analyses the sequence does not know the real content
of the box.) 
} 
Thus after 100 Black in a row 
we become 
`practically certain' to observe Black
in the 101-th extraction, being the probability of White only 
$4\times 10^{-11}$ -- but yet
{\em not impossible!}\,\footnote{You can 
evaluate this probability also  using the 
 of R code in  Footnote 31 of Ref. \cite{ME2016} remember to
change the first \verb|sprintf()| format, for example with 
\verb|%.15f| or \verb|%.5e|. The results are
`numerically' the same.}
(If you think that very improbable events do not occur
in real life, then wait for Appendix B.)

For comparison we could show in the same plot also the relative
frequency of White in the $n$ extractions, as we shall do with the
simulations from the other boxes, but, since in this case
it is always zero, it is of little interest, and 
anyway not representable in a log scale. 
It is, instead, more interesting, the probability
evaluated
({\em incorrectly}!) applying the 
Laplace rule of succession (Equation 15 of Ref. \cite{ME2016}),
that in this case becomes $1/(n+2)$, and, 
by complement $(n+1)/(n+2)$ for Black. 
As we can see, the performance is rather poor.

However is not Laplace
to be wrong\footnote{This is, for example, verbatim what he wrote
concerning his too much misunderstood probability
of the sun rising tomorrow: ``{\sl On trouve ainsi qu'un 
\'ev\'enement \'etant arriv`'e de suite, un nombre quelconque de fois;
la probabilit\'e qu'il arrivera encore la fois 
suivente, est \'egale \`a ce nombre augment\'e de l'unit\'e, 
divis\'e par le m\^eme nombre augment\'e
de deux unit\'es.} [$(n+1)/(n+2)$] {\sl En faisant, 
par example, remonter la plues ancienne \'epoque de l'histoire,
\`a cinq mill ans, ou \`a 1826213 jours, et 
le solei s'\'etant lev\'e constenmment dans
cet intervall, \`a chaque r\'evolution de vinght-quatre heures; il y a
1826214 \`a parier contre un, qu'il se levera encore demain.
{\it Mais ce nombre est incomparablement plus fort pour celui
qui connaissant par l'ensemble des ph\'enom\`enes, 
le principe r\'egulateur des jours et des saisons, 
voit que rien \underline{dans le moment actuel}, 
ne peut en arr\^eter le cours.}}'' [{\sl ``Thus we find that 
an event having occurred successively
any number of times, the probability that it will
happen again the next time is equal
to this number increased by unity divided
by the same number, increased by two units. Placing the most ancient epoch
of history at five thousand years ago, or 
at 182623 days, and the sun having risen constantly in the interval
of each revolution of twenty-four hours, it is a bet of 1826214 
to one that it will rise again tomorrow.} {\it But this number
 is incomparably greater for him who, recognizing in the totality
 of phenomena the principal regulator of days and seasons, 
sees that nothing \underline{at the present moment}
 can arrest the course of it}.''] 
(italics and underline mine) \cite{Laplace}. Great Laplace! 
(And please note once more the probability
expressed in terms of a virtual coherent bet.)
\label{fn:Laplace_Sole} 
}, 
but rather those who would use his formula
a-critically, without understanding the assumptions
behind it, which were discussed in detail in the text. 
In our specific case, as it might be in important cases 
of real life, the prior of the propensity of the box to give
White was not uniform between 
between 0 and 1. We had instead only six possible values,
and the full calculation takes into account of the real 
situation.\footnote{I have called the attention several times
in the past (e.g. \cite{BR}) that prior-less methods
are not by default `objective' because ``they do not use
priors.'' On the contrary, it is possible to show that 
there are in most cases hidden (most times) uniform 
priors, like in the result of the so called maximum-likelihood 
method of the statisticians (see e.g. \cite{BR}).
} For this reason the name of Laplace is in quote marks
in the legends of the figures, to mean ``misused Laplace rule.''

\section{Box $B_1$}
\begin{figure}
\centerline{\includegraphics[width=0.9\linewidth]{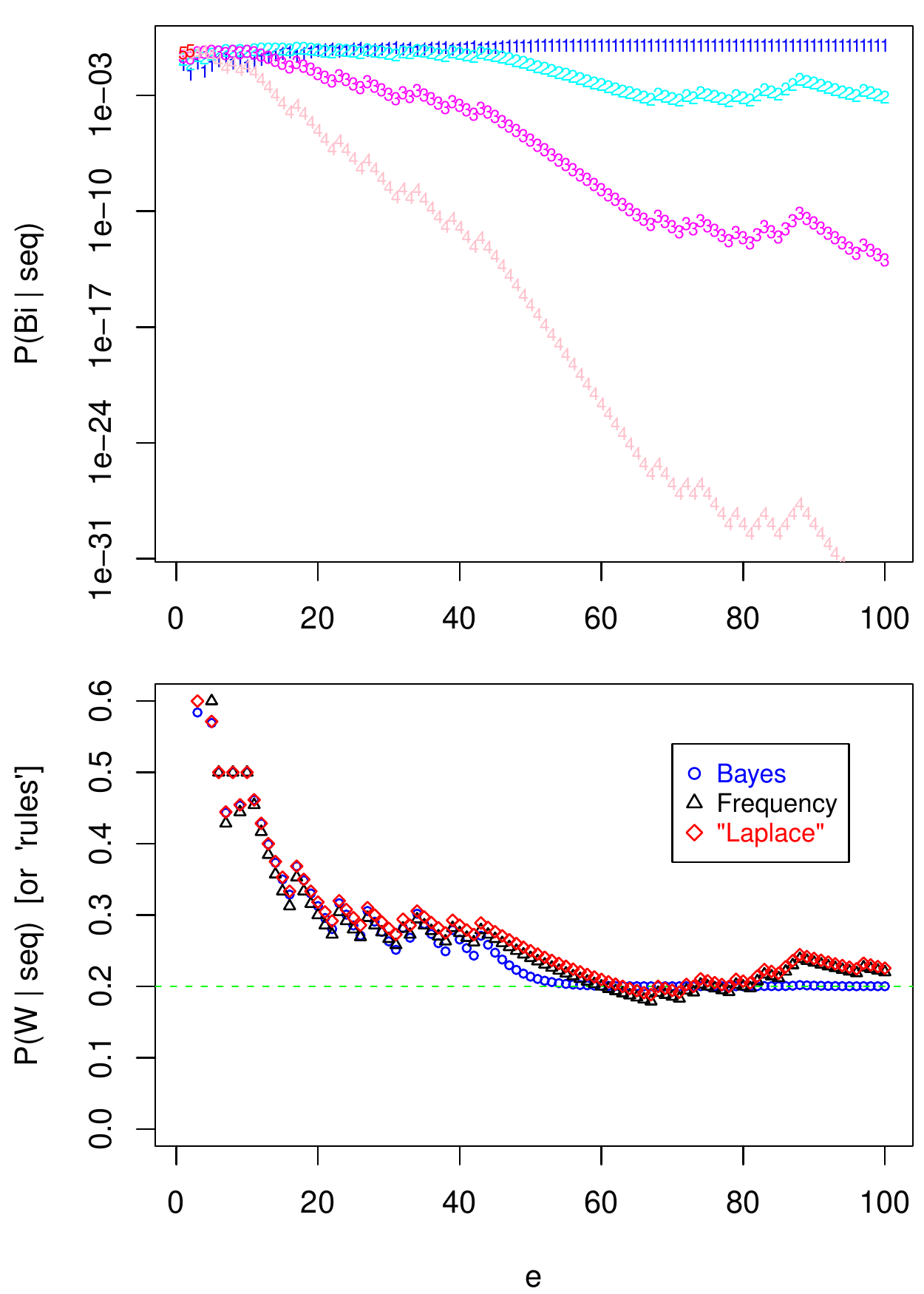}}
\caption{\small \sf Simulated extractions from $B_1$ (run 1: 1:100).}
\label{fig:Simulazione_B1_0_100}
\end{figure}

\begin{figure}
\centerline{\includegraphics[width=0.9\linewidth]{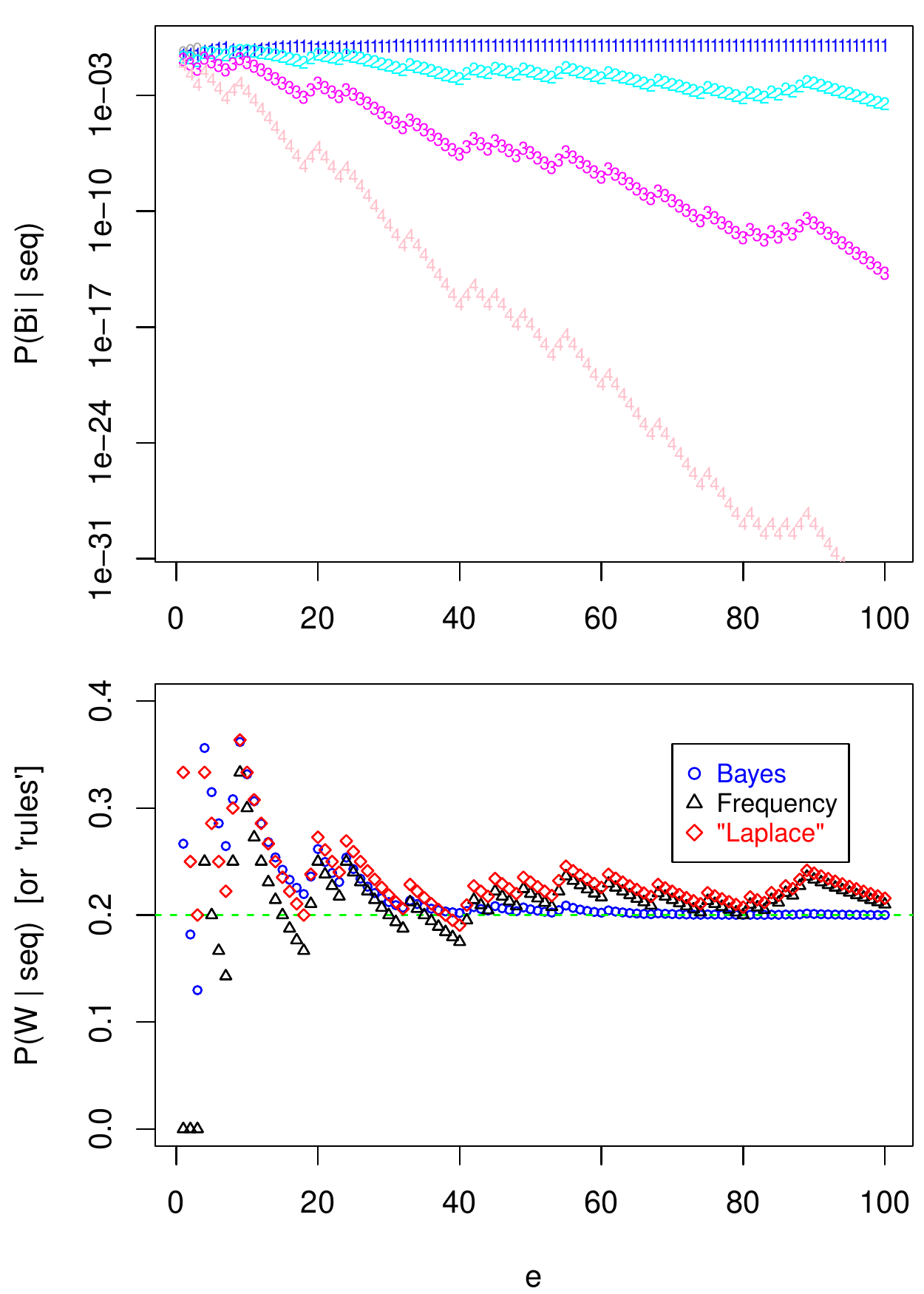}}
\caption{\small \sf Simulated extractions from $B_1$ (run 2: 101:200).}
\label{fig:Simulazione_B1_100_100}
\end{figure}

\begin{figure}
\centerline{\includegraphics[width=0.9\linewidth]{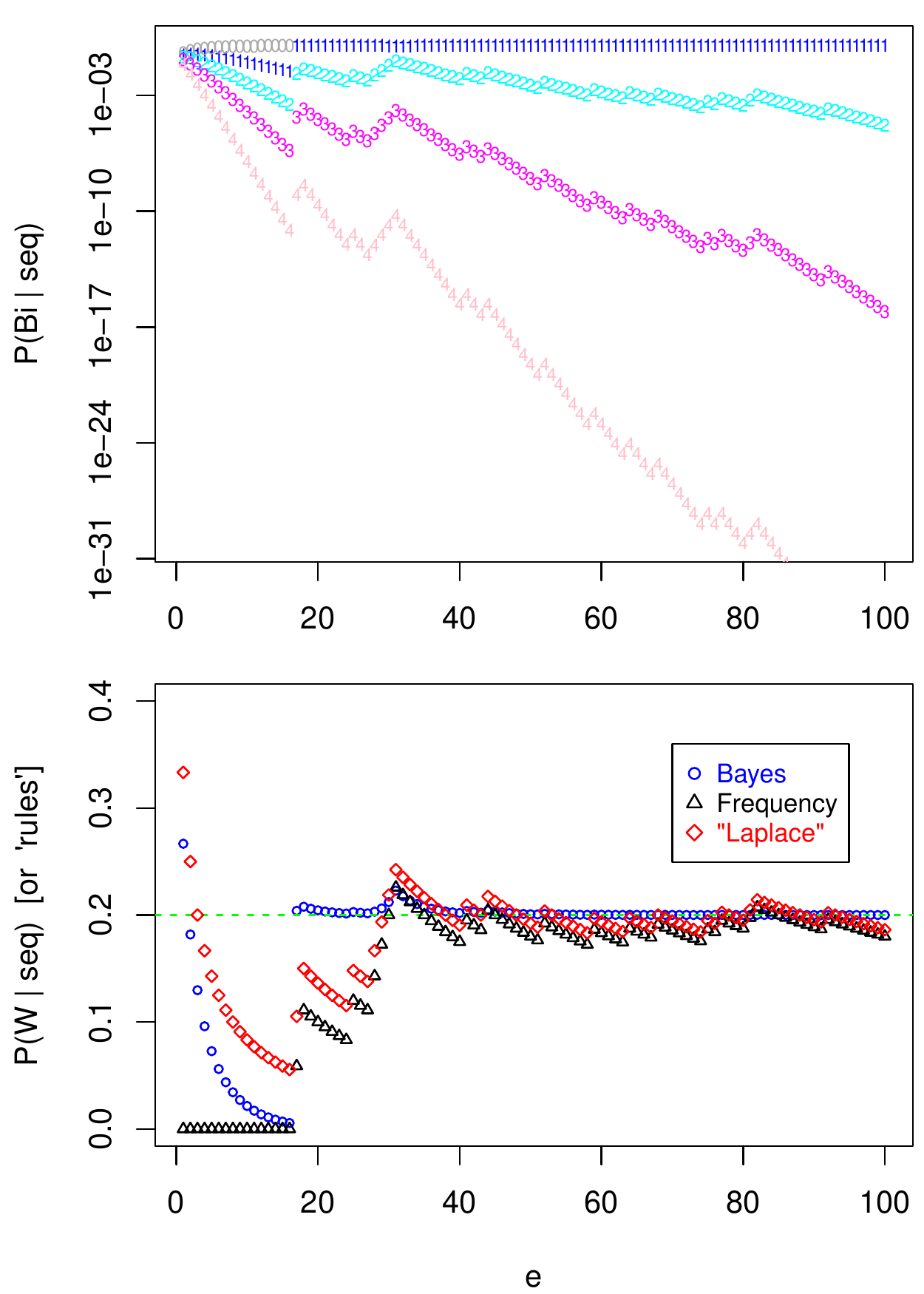}}
\caption{\small \sf Simulated extractions from $B_1$ (run 3: 201:300).}
\label{fig:Simulazione_B1_200_100}
\end{figure}

The analysis of the sequences 
from box $B_1$ (and, by symmetry, from $B_4$) 
is in general the most interesting and instructive, 
because the probabilities calculated using 
probability theory, taking into account all the available information
properly, differ quite a lot from those obtained 
using intuitive heuristics, or  from `prescriptions'
based on improper use of theoretical results not fully understood
(see Footnote \ref{fn:Laplace_Sole}). 
Figures 
\ref{fig:Simulazione_B1_0_100}, 
\ref{fig:Simulazione_B1_100_100}
and \ref{fig:Simulazione_B1_200_100} show the results
of the inferences and of the (probabilistic) predictions
based on three sequential runs of 100 extractions each. 
Each story is peculiar, 
as real life situations are, 
and we see that -- in the simulations 
we know the `truth' -- the method based
on probability theory, and which take into account at best 
all available information,
performs {\em much better} than the others. 

It is worth 
remembering that all real cases
are {\em unique} and we can only rely on the quality
of the methods, as well as of the data and all relevant information.
As someone says (reference missing), in the Bayesian 
analysis ``the result is the result.'' For example, in the
first part of  
sequence on which of Fig.\,\ref{fig:Simulazione_B1_0_100}
is based, $B_1$ and  $B_2$ seemed practically equally likely, 
and, as consequence, the probability of White in the 
next extraction was in between 1/5 and 2/5. That's all. 
This the best we could say at that moment, but as soon as the 
overall relative frequencies of White approaches 20\%
(frequencies are reported as black triangles in the lower plot)
there as a kind of `attraction' from $B_1$: its probability
suddenly rises, and the probability of White approaches rapidly
20\%. Once balls of both colors are observed,
if the relative frequency of observed White 
goes under 20\% the effect of `attraction'
gets more enhanced, because the next possibility, related
to box $B_2$, is too far.\footnote{No esoteric meaning is attached
  to the term `attraction'. It is just because
  the next possible value of propensity, 4/5 of box $B_2$, is ``too far''
  -- see Appendix B.
} 
  
Also interesting is run 3 (Fig.\,\ref{fig:Simulazione_B1_200_100}),
characterized by 16 black in a row.\footnote{But if you check the file
you will see that there were already 11 Black just before, 
summing thus to 27 Black in a run (and, after 2 White, 
other 6 Black follow). It simply happened, and 
for this reason I would like to insist on the worries already
expressed in Footnote 28 of Ref. \cite{ME2016}, i.e. interpreting 
probabilistic statements as pedantic regularities. 
} As a result, for a while we believed stronger and stronger
we had picked up $B_0$, and thus the probability of White
in the next extraction was (exponentially) decreasing. 
resulting in small probability
of White in the next extraction. Than, suddenly,  
we observe White, and  the probability of $B_0$ instantly drops 
to zero,\footnote{I am sorry for those who dislike discontinuities,
but a crucial single (solid) experimental information 
can change dramatically our vision of the world, as it happens
to those who suddenly learn that their quite and polite neighbor was 
indeed a serial killer keeping rests of human bodies in his fridge.} 
while the probability $B_1$ jumps practically at 100\% 
and that of a next White at 20\%.\footnote{More precisely,\\
\verb|> N=5; i=0:N; pii=i/N; n=17; x=1|\\
\verb|> ( PBi =  pii^x * (1-pii)^(n-x) / sum( pii^x * (1-pii)^(n-x) ) )|\\
\verb|[1] 0.000000e+00 9.803047e-01 1.965040e-02 4.487479e-05 9.129799e-10|\\
\verb|> sum( pii * PBi )|\\
\verb|[1] 0.203948|
}
It is nice, and instructive, to observe that from this extraction
on, $P(\mbox{W}\,|\,\mbox{seq},I)$ will always be 
{\em slightly above 20\%}. Those who have a biased mind would
 speak about a `biased estimator'. In reality, it is 
a just logical consequence of the fact that, once we have ruled out $B_0$, 
the probability of White in a future extraction, 
weighted average of all possible values of the propensity
of the box to give White,  {\em has to be} slightly 
above the minimum possible value of propensity, that  
is 1/5. The frequency based value, as well as that from  the
{\em \underline{misused} Laplace rule}, remains quite for a while 
below 20\%, and than oscillates around it, in contradiction
with the fact that $B_0$ has been definitively ruled out.  

\begin{figure}
\centerline{\includegraphics[width=0.9\linewidth]{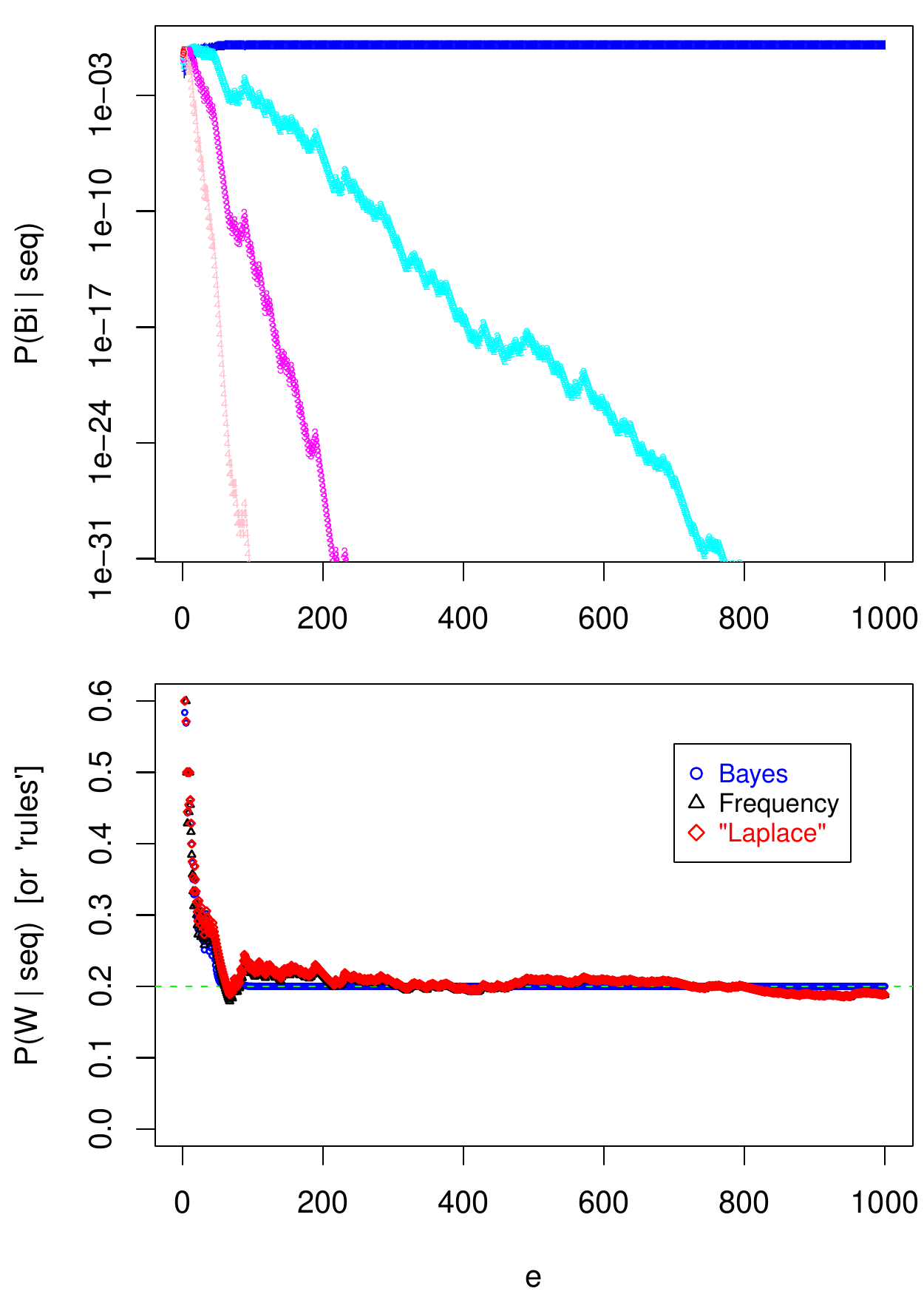}}
\caption{\small \sf Simulated extractions from $B_1$ (full sequence 1:1000).}
\label{fig:Simulazione_B1_0_1000}
\end{figure}

You might object that after a very long sequence 
also the other evaluations will eventually 'converge' 
(see Footnote 28 of Ref. \cite{ME2016} 
for remarks on the precise meaning of this
term in probability theory), but, as it as been famously said, 
``{\sl in the long run we are all 
dead}.''\footnote{\url{https://en.wikiquote.org/wiki/John_Maynard_Keynes#Quotes}} 
Let us see what happens if we 
analyze the full sequence of 1000 extractions 
(Fig. \ref{fig:Simulazione_B1_0_1000}). 
The frequency based evaluations of the next 
observation
is still oscillating around 20\%, while that obtained from 
probability theory approaches 1/5 (from above!)
by $10^{\approx -40}$ [rough estimate obtained extrapolating the probability
of $P(B_2\,|\,\mbox{seq})$ from the above plot]. 

\section{Box $B_2$}
\begin{figure}
\centerline{\includegraphics[width=0.9\linewidth]{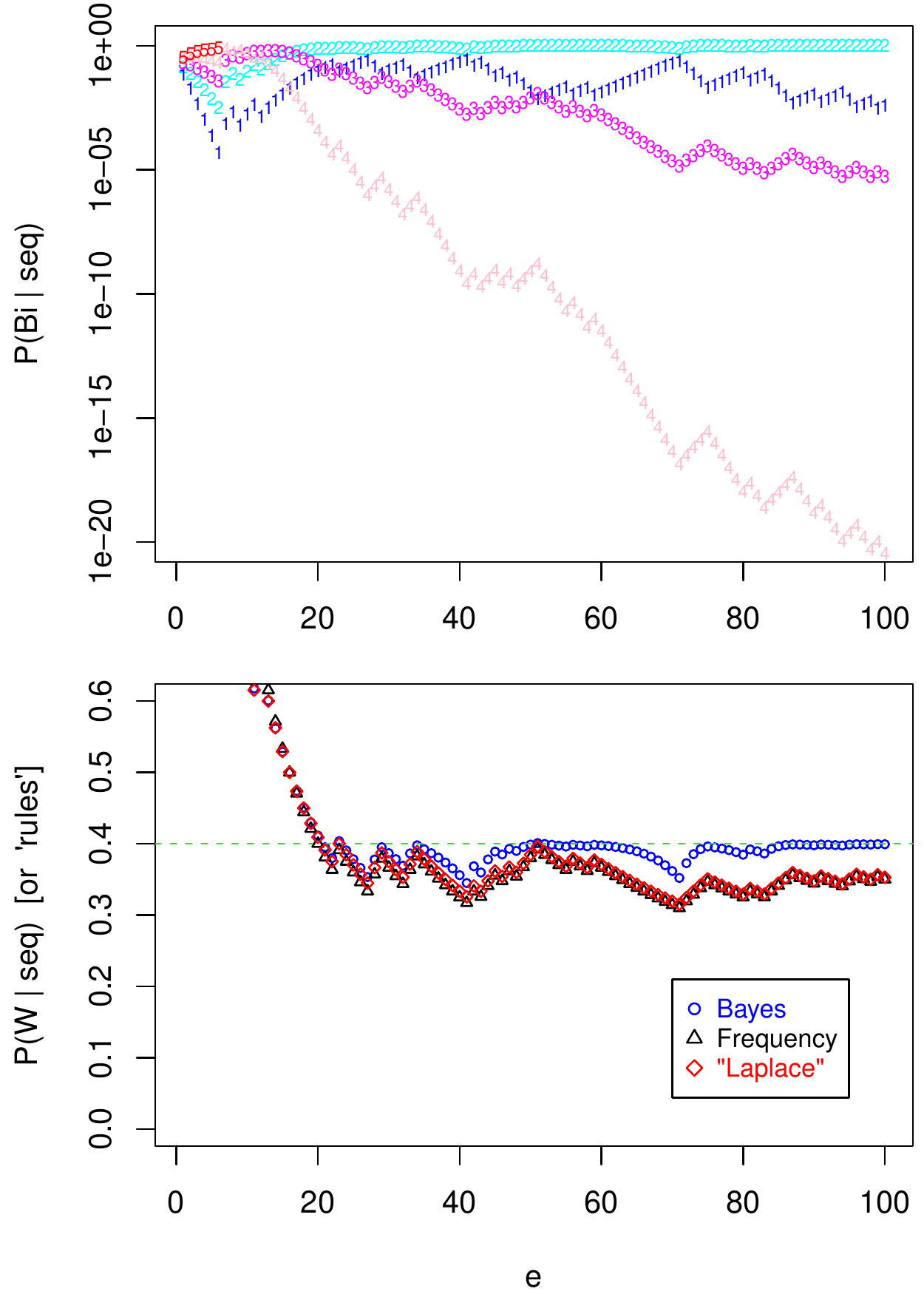}}
\caption{\small \sf Simulated extractions from $B_2$ (run 1: 1:100).}
\label{fig:Simulazione_B2_0_100}
\end{figure}

\begin{figure}
\centerline{\includegraphics[width=0.9\linewidth]{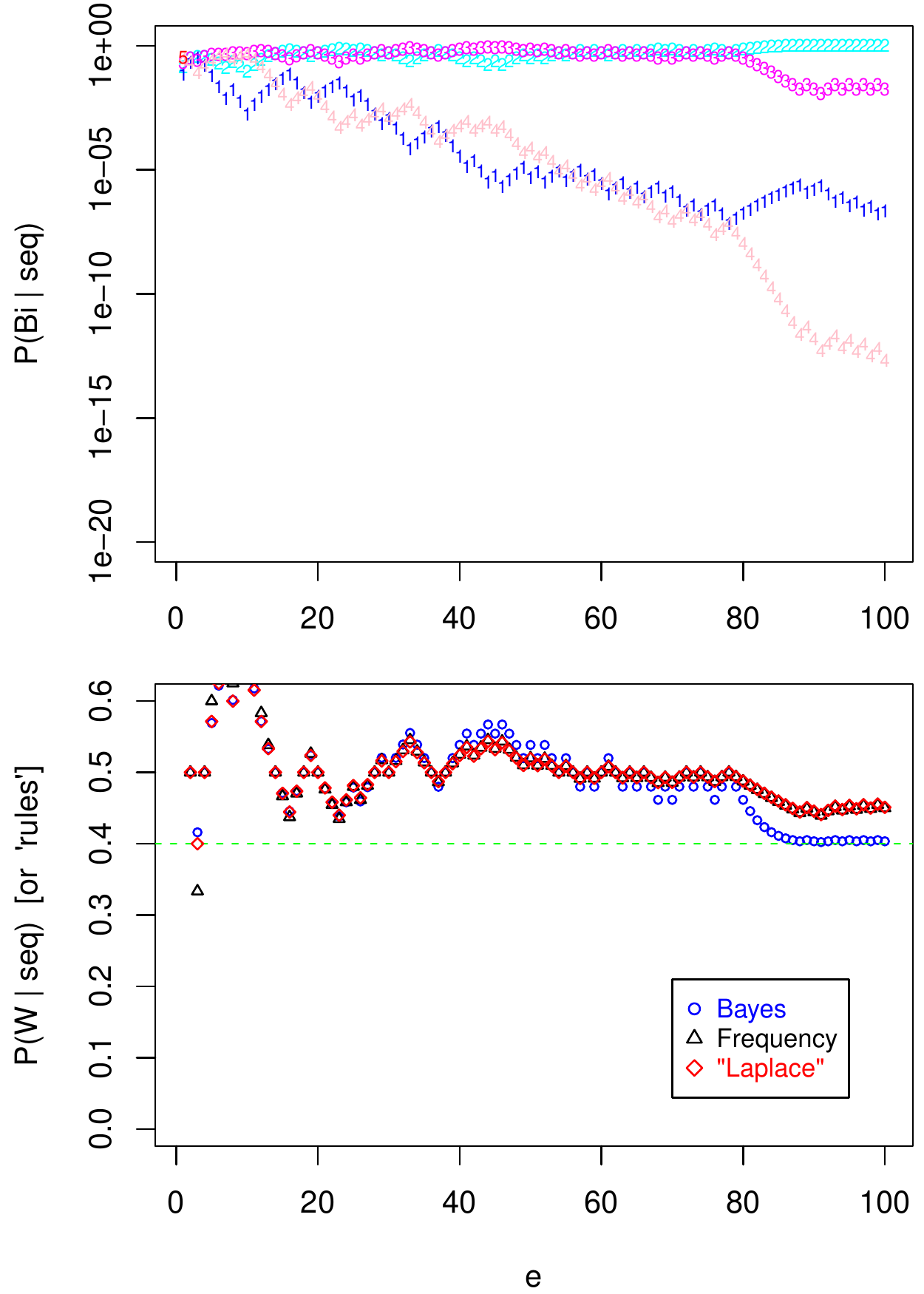}}
\caption{\small \sf Simulated extractions from $B_1$ (run 2: 101:200).}
\label{fig:Simulazione_B2_100_100}
\end{figure}

\begin{figure}
\centerline{\includegraphics[width=0.9\linewidth]{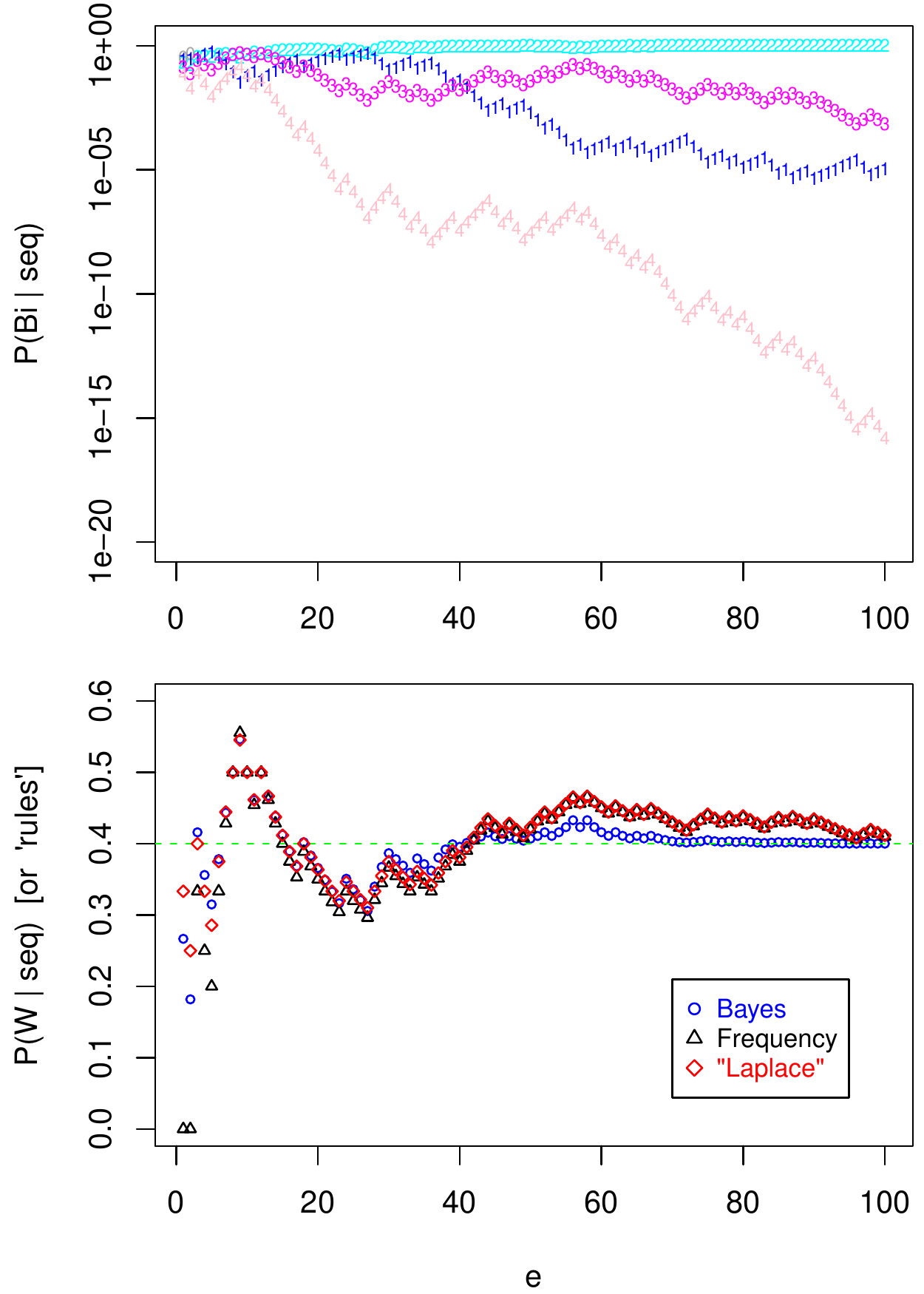}}
\caption{\small \sf Simulated extractions from $B_2$ (run 3: 201:300).}
\label{fig:Simulazione_B2_200_100}
\end{figure}

Let us conclude this round up by also showing 
the results of the analysis of three
runs and of the complete sequence from the box $B_2$ 
(Figures \ref{fig:Simulazione_B2_0_100}, 
 \ref{fig:Simulazione_B2_100_100},  
\ref{fig:Simulazione_B2_200_100} and  \ref{fig:Simulazione_B2_0_100})
\begin{figure}
\centerline{\includegraphics[width=0.9\linewidth]{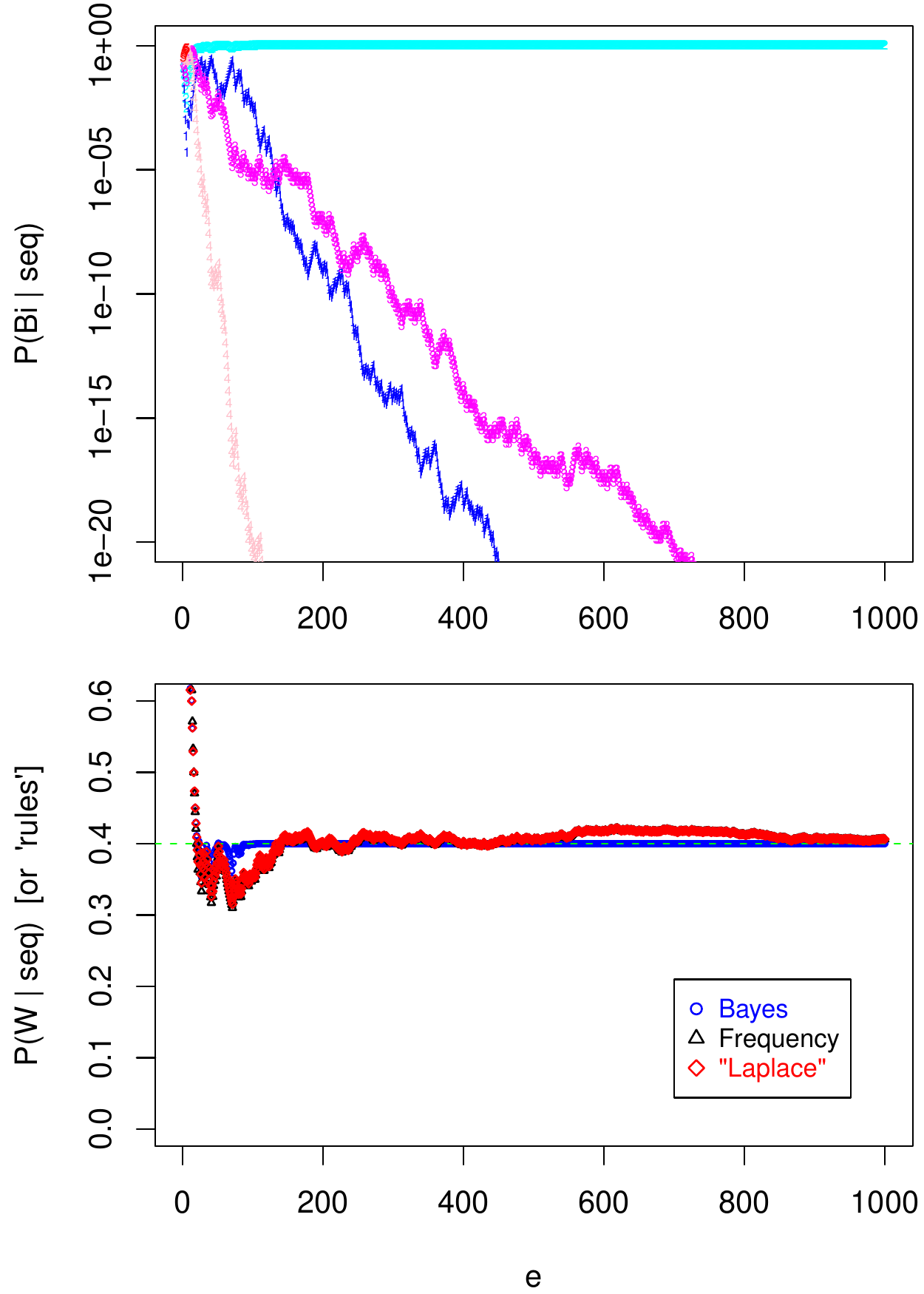}}
\caption{\small \sf Simulated extractions from $B_2$ (full sequence 1:1000).}
\label{fig:Simulazione_B2_0_1000}
\end{figure}
without further comments, besides that in probabilistic
matter, as in real life, 
the German dictum ``einmal ist keinmal'' applies. 

\section{$P(B_i\,|\,\mbox{seq},I)$: Bayesian vs frequentistic 
evaluations  }
After having seen the comparisons of the evaluations of probabilities of 
White in a future extraction, someone would like to see
something similar concerning the probability of the different
box compositions in the light of the observed sequence. 
But he/she will be disappointed to learn that such a comparison
is simply impossible because {\em frequentists prohibit the
very concept of} $P(B_i\,|\,\mbox{seq},I)$. That's all, 
Sorry! (And I am sorry for you, if you thought you were a frequentist,
but, nevertheless, you also thought that
such a probability had a 
meaning -- see \cite{BadMath, WavesSigmas} and references therein
for details).\\  
\mbox{} \\

\section{Most of our observations {\em had} 
very small chance to occur}
An important misconception about probability is to
confuse the probability of the effect given a given 
hypothesis with the probability of that hypothesis
given the observed effect. The name ``prosecutor fallacy'',
with which this logical error is often designed 
gives, by itself, an idea of its relevant importance 
in real life.\footnote{See  e.g.
\url{http://www.agenarisk.com/resources/probability_puzzles/prosecutor.shtml}.
This fallacy is somehow similar to the misinterpretation
of `p-values' as probability of hypotheses 
(see \cite{BadMath,WavesSigmas} and references therein), 
but, even if the numbers are less impressive, 
the logical fallacy of misinterpreting p-values
 is even worst, because the erroneous 
conclusion is not based solely on the data, 
but also on data less probable than those actually observed.
}
Continuing with the style of this paper, I would like 
to touch this point using the third run of the sequence
from box $B_1$ (Fig. \ref{fig:Simulazione_B1_200_100}),
which I find particular instructive. We shall analyze what we have learned 
after  the 16th, the 17th and the 100th extraction, 
also giving the details of the calculations in R,
which start with the usual initialization 
  (\verb|N=5; i=0:N; pii=i/N|):
\begin{description}
\item{$n=16$ (run 3):} At the beginning we had 16 black in a row, 
              resulting on the following probabilities:\\
{\small 
              \verb|> n=16; x=0; ( PBi =  pii^x * (1-pii)^(n-x) / sum( pii^x * (1-pii)^(n-x) ) )|\\
              \verb|[1] 9.723559e-01 2.736939e-02 2.743123e-04 4.176237e-07 6.372432e-12|\\
              \verb|> sum( pii * PBi )|\\
              \verb|[1] 0.005583852| \\
}
              We are 97\% confident to have got $B_0$, 
              2.7\% to have got $B_1$, and so on. On the other hand,
              the probabilities to get ``0 white in 16 trials'' -- be careful,
              I am trying to fool you -- under the different hypotheses $B_i$
              are $((5-i)/5)^{16}$, that is\\
              {\small 
              \verb|> ((5-i)/5)^16|\\
              \verb|[1] 1.000000e+00 2.814750e-02 2.821110e-04 4.294967e-07 6.553600e-12 0.000000e+00| 
              }\\ 
              So it seems than that the small probability to $B_1$ is 
              due to the small probability to get the `observation' 
              $(x=0,n=16)$ from $B_1$; and similarly with the other
              boxes which contain white balls. 
\item{$n=17$ (run 3):}
              Here is what happened after the next extraction, in which
              we observe White:\\
              {\small 
              \verb|> n=17; x=1; ( PBi =  pii^x * (1-pii)^(n-x) / sum( pii^x * (1-pii)^(n-x) ) )|\\
              \verb|[1] 0.000000e+00 9.803047e-01 1.965040e-02 4.487479e-05 9.129799e-10  0.000000e+00|\\
              \verb|> sum( pii * PBi )|\\
              \verb|[1] 0.203948|
              }\\
              What is, instead, the probability of the observation, subject
              to the different compositions? 
              You might think at binomial distributions resulting in 1 success
              in 17 trials, with probabilities given by $i/5$, that is\\
              {\small 
              \verb|> dbinom(x, n, pii)|\\
              \verb|[1] 0.000000e+00 9.570149e-02 1.918355e-03 4.380867e-06 8.912896e-11 0.000000e+00|
              }\\
              However, this is not we have really observed, but just 
              {\em its summary}. The observation was indeed {\em the} sequence!
              And the probability of the sequence is quite 
              different:\footnote{The difference is due to the binomial 
               coefficient, 
               equal to 17 in this case,
               for which we have to divide the previous numbers.}\\
              {\small 
              \verb|> pii^x * (1-pii)^(n-x)|\\
              \verb|[1] 0.000000e+00 5.629500e-03 1.128444e-04 2.576980e-07 5.242880e-12 0.000000e+00|
              }
\item{$n=100$ (run 3):} As it easy to predict, the difference becomes `dramatic' for very large values of $n$.
               Having observed 18 White in 100 extractions, 
               these are the probabilities of the hypotheses:\\
              {\small 
              \verb|> n=100; x=18; ( PBi =  pii^x * (1-pii)^(n-x) / sum( pii^x * (1-pii)^(n-x) ) )|\\
              \verb|[1] 0.000000e+00 9.999851e-01 1.491273e-05 8.011548e-17 2.938692e-39 0.000000e+00|
              }\\
              to be compared with the probability of 18 successes in 100 trials 
              for the different boxes:\\
              {\small 
              \verb|> dbinom(x, n, pii)|\\
              \verb|[1] 0.000000e+00 9.089812e-02 1.355559e-06 7.282455e-18 2.671256e-40|
              }\\
              But the conditional probabilities of what we have really observed are now strikingly different:\\
              {\small 
              \verb|> pii^x * (1-pii)^(n-x)|\\
              \verb|[1] 0.000000e+00 2.964277e-21 4.420612e-26 2.374881e-37 8.711229e-60  0.000000e+00|
              }\\
              But indeed, the fact that this sequence {\em had}  $3\times 10^{-21}$ {\em chance} 
              (really in the sense of a propensity)
              to occur from 
              $B_1$ {\em is absolutely irrelevant}. 
              What matters is only this probability
              with respect to those from the other boxes.
              Indeed the respective conditional probabilities 
              provide the Bayes-Turing factors for every pair of hypotheses. 
              And therefore, since {\em in our toy experiment} the different compositions were initially
              {\em equally likely}, we get odds of $B_1$ vs $B_2$ of $6.7\times 10^4$;
              vs  $B_3$ and  $B_4$  $1.2\times 10^{16}$ and  $3.4\times 10^{38}$, respectively. These are the numbers that matter.
\end{description}
At this point I hope the lesson is clear, without you need  to be further
impressed with the numbers that we would get analyzing the full sequence 
of 1000 extractions:
\begin{itemize}
\item that fact that we can make our inference and prediction based on the 
      number of trials and the number of successes it is because 
      these {\em summaries} are {\em `sufficient'} 
      for the purpose of the inference (and forecasting);
      but the real {\em observation}
      is the sequence;
\item most of the fact of real life {\em had} very little chance to occur,
      if we analyze them with enough accuracy. But this implies little on the 
      probabilities of the cause that might have produced them. 
      What matters are the ratio of conditional probabilities:
      $P(E\,|\,C_i,I)/P(E\,|\,C_j,I)$.
\end{itemize}

\section{Conclusions}
Having to evaluate probabilities of hypotheses
and probabilities of occurrences of future events,
unless you possess a crystal ball, 
it is hard to out-perform Bayesian reasoning,
if it is used consistently, and all the available
pieces of information are properly taken into account.
But the lesson which comes from playing with the simulated
sequences goes beyond the demonstration of the power of the so called
Bayesian methods.

For example
I find it important that, in the
training of {\em probabilistic thinking},
people should be exposed to the rich 
variety of what can occur 
randomly. And, therefore, most events
of real life {\em had} very little chance
to occur. Think, for example, at a given configuration 
of light content in pixels, when 
you shoot a picture with a digital camera. More simply,
and easier to calculate, consider a number to twelve decimal places
that can come from a Gaussian random number generator, 
like that obtained with the following R commands:\\
{\small 
\verb|> options(digits=14); set.seed(20160715); nd=12; dxm=10^(-nd)/2|\\
\verb|> (xr=round(rnorm(1), nd)); as.double(sprintf("%.2e", dnorm(xr,) * 2 * dxm))|\\
\verb|[1] 1.479427401471|\\
\verb|[1] 1.34e-13|\\
\verb|> (xr=round(rnorm(1), nd)); as.double(sprintf("%.2e", dnorm(xr,) * 2 * dxm))|\\
\verb|[1] -0.762658301757|\\
\verb|[1] 2.98e-13|
}\\
(Yes, every time you repeat this line of code you will observe, 
{\em with certainty}, numbers which {\em had} about 1-in-trillions
chance to occur! And they all come with probability 1
from a Gaussian random generator with $\mu=0$ and and $\sigma=1$)

The reason I insist on these apparently trivial considerations
is that I have seen too often in the past, and even quite recently, 
attempts of indoctrinating people with `statistical regularities'.
These attempts imply a misinterpretation of probability theorems and, 
at the same time, a refusal of the concept 
of a single event probability.
Instead, not only degrees of beliefs apply to single events,
but also propensities, if we reflect on the fact that
a propensity might change with time \cite{ME2016}.

\section*{Acknowledgements}
This work was partially supported by a grant from Simons Foundation
which allowed me a stimulating working environment during
 my visit at the Isaac Newton Institute
of Cambridge, UK (EPSRC grant EP/K032208/1). It is a pleasure
to thank Paolo Agnoli of Pangea Formazione, 
for discussions on probabilistic matter, 
always with some philosophical flavor, encouragements to go 
ahead with the six box toy experiment, 
and for comments  on the manuscript, which 
has also benefitted of comments by Patricia Wiltshire.

\end{document}